\documentstyle[12pt]{article}

\def\be{\begin{equation}}
\def\ee{\end{equation}}
\def\bea{\begin{eqnarray}}
\def\eea{\end{eqnarray}}

\def\1{\'{\i}}
\def\R{{\rm I\kern-.2em R}}

\def\picture #1 by #2 (#3){
  \vbox to #2{
    \hrule width #1 height 0pt depth 0pt
    \vfill
    \special{picture #3} 
    }
  }

\def\scaledpicture #1 by #2 (#3 scaled #4){{
  \dimen0=#1 \dimen1=#2
  \divide\dimen0 by 1000 \multiply\dimen0 by #4
  \divide\dimen1 by 1000 \multiply\dimen1 by #4
  \picture \dimen0 by \dimen1 (#3 scaled #4)}
  }

\begin{document}

\thispagestyle{empty}

\hfill \ \quad
\
\vspace{2cm}

\begin{center} {\LARGE{\bf{Universal $R$-matrix for }}}

 {\LARGE{\bf{non-standard quantum $sl(2,\R)$}}}
 \end{center}

\bigskip\bigskip\bigskip

\begin{center} Angel Ballesteros and  Francisco J. Herranz
\end{center}

\begin{center} {\it {  Departamento de F\1sica, Universidad
de Burgos} \\   Pza. Misael Ba\~nuelos  \\
E-09001, Burgos, Spain}
\end{center}

\bigskip\bigskip\bigskip

\begin{abstract}
A universal $R$-matrix for the non-standard (Jordanian)  quantum
deformation of
$sl(2,\R)$ is presented. A family of solutions of
the quantum Yang--Baxter equation is obtained from
some finite dimensional representations of this Lie  bialgebra
quantization of
$sl(2,\R)$.
\end{abstract}

\newpage

The quantum Yang--Baxter equation (YBE)
\be
{\cal  R}_{12}{\cal  R}_{13}{\cal  R}_{23}=
{\cal  R}_{23}{\cal  R}_{13}{\cal  R}_{12}.
\label{aac}
\ee
was discovered to play a relevant role as the integrability
condition for (1+1) quantum field theories \cite{Yang}  and also
in connection with two dimensional models in lattice statistical
physics \cite{Baxter}, conformal field theory \cite{Gomez} and
knot theory \cite{Turaev}. Nowadays it is well known that the
investigation on the algebraic properties of this equation and
the obtention of new solutions are closely related to the study
of quantum groups and algebras \cite{FRT}.

In fact, let ${\cal A}$ be a Hopf algebra and let  ${\cal R}$ be
an invertible element in ${\cal A}\otimes {\cal A}$ such that
\be
\sigma \circ \Delta( X) ={ \cal  R} \Delta( X)
{\cal  R}^{-1},  \qquad
\forall
\, X \in {\cal  A} \label{aaa}
\ee
where $\sigma$ is the
flip operator $\sigma(x\otimes y)=(y\otimes x)$. If  we write
${\cal  R}=\sum_{i} a_i\otimes b_i$, ${\cal
R}_{12}\equiv\sum_{i} a_i\otimes b_i\otimes 1$, ${\cal
R}_{13}\equiv\sum_{i} a_i\otimes 1\otimes b_i$,
${\cal  R}_{23}\equiv\sum_{i} 1\otimes a_i\otimes b_i$  and the
relations
\be
(\Delta  \otimes id){\cal  R}
={\cal  R}_{13}{\cal  R}_{23},\qquad
  (id \otimes \Delta){\cal  R} ={\cal  R}_{13}{\cal  R}_{12},
\label{aab}
\ee
are fulfilled, $({\cal  A},{\cal  R})$ is called  a {\em
quasitriangular} Hopf algebra \cite{Dr}. In that case, ${\cal
R}$ is easily proven to be a solution of (\ref{aac}). Hereafter,
if ${\cal A}$ is a Hopf algebra and ${\cal  R}$ fulfills both
(\ref{aac}) and (\ref{aaa}), we shall say that
${\cal  R}$ is a (quantum) {\em universal $R$-matrix}  for ${\cal
A}$. Obviously, different representations for the algebra ${\cal
A}$ will give rise to different explicit solutions of the quantum
YBE.

In particular, let ${\cal A}$ be a quantum  deformation $U_z (g)$
of the universal enveloping algebra of a Lie algebra $g$. Then,
the quasicocommutativity property (\ref{aaa}) is translated, in
terms of the Hopf algebra dual to ${\cal A}$, into the known FRT
relations defining the quantum group $Fun_z(G)$ \cite{FRT}. It is
also known that each $U_z (g)$ defines a unique Lie bialgebra
structure on
$g$ that can be used to characterize the quantum deformation. For all
semisimple Lie algebras, all these Lie bialgebra structures  are
coboundaries generated by  classical $r$-matrices. In $sl(2,\R)$,
two outstanding Lie bialgebra structures can be mentioned: the
standard one, generated by the classical $r$-matrix $r^{(s)}=
\lambda\,J_+\wedge J_-$
 and the non-standard (triangular) Lie bialgebra given
by the element  $r^{(n)}=\chi\,J_3\wedge J_+$
($\lambda,\chi\,\in\,\R$). The
quantization of the former is the well known  Drinfel'd--Jimbo
deformation, whose quantum universal $R$-matrix was given in
\cite{Kirillov}. For the latter, the corresponding non-standard
quantum algebra is the so called Jordanian deformation of
$sl(2,\R)$ \cite{Ohn} (introduced for the first time in
\cite{Demidov} in a quantum group setting; see  \cite{AKS} and
references therein and \cite{Beyond}, where the non-standard
deformation of $so(2,2)$ was constructed). To our knowledge, no
quantum universal
$R$-matrix is known for this case  (the $R$ given in \cite{Ohn}
is neither a solution of (\ref{aac}) nor verifies (\ref{aaa}), as
it has already been pointed out in
\cite{AA,RR}). The aim of this letter is to provide it.

Let us consider the $sl(2,\R)$  Lie  algebra with the following
commutation relations
\be
[J_3,J_+]=2J_+,\qquad  [J_3,J_-]=-2J_-,\qquad [J_+,J_-]=J_3 .
\label{aa}
\ee
The classical $r$-matrix
\be
r= z\, J_3\wedge J_+
\label{ab}
\ee
is a solution of the classical YBE and it generates the Lie
bialgebra structure with cocommutators given, as usual, by
$\delta(X)=[1\otimes X+X\otimes 1,r]$:
\be
\delta(J_+)=0,\qquad \delta(J_3)=2\,z\, J_3\wedge J_+,\qquad
\delta(J_-)=2\,z\, J_-\wedge J_+ .
\label{ac}
\ee
This map is the first order (in the
deformation parameter $z$) of the co-antisymmetric  part of the
deformed coproduct corresponding to the non-standard quantization
of $sl(2,\R)$    given by the relations \cite{Ohn}
\bea
&&  \Delta (J_+)  =1 \otimes J_+  + J_+ \otimes 1,\cr
&& \Delta (J_3) =e^{-z J_+ } \otimes J_3 + J_3\otimes
e^{z J_+ }, \label{ad} \\
&& \Delta (J_-) =e^{-z J_+ } \otimes J_- + J_-\otimes
e^{z J_+ },
\nonumber
\eea
\bea
&& [J_3,J_+ ]=2\,\frac{\sinh ( z J_+  )} z,\qquad
[J_+  ,J_- ]= J_3,\\
&& [J_3,J_-]=-J_-\cosh (z J_+  ) -\cosh (z J_+  ) J_- .\quad
\label{ae}
\eea

In order to
obtain a universal $R$-matrix linked to this quantum  algebra, an
essential point is to choose an adequate basis for the
deformation. Let us perform a (non linear) transformation of the
generators $\{J_3,J_+,J_-\}$ as follows:
\bea
&&A_+=J_+,\qquad A= e^{z J_+ }J_3,\cr
&&A_-= e^{z J_+ }J_- - \frac z4 e^{z J_+ } \sinh( z J_+  ).
\label{af}
\eea
This change of basis leads to a Hopf algebra  (here denoted by
$U_z sl(2,\R)$) with the following  coproduct
$(\Delta)$, counit $(\epsilon)$, antipode  $(\gamma)$ and
commutation rules:
\bea
&&  \Delta (A_+ ) =1 \otimes A_+  + A_+ \otimes 1,\cr
&& \Delta (A) =1 \otimes A + A\otimes
e^{2z A_+ }, \label{ag} \\
&& \Delta (A_-) = 1 \otimes A_- + A_-\otimes
e^{2 z A_+ },
\nonumber
\eea
\be
\epsilon(X) =0,\qquad
 \mbox{for $X\in \{A,A_+,A_-\}$},
 \ee
\be
\gamma(A_+)=-A_+,\qquad \gamma(A)=-A e^{-2zA_+},\qquad
\gamma(A_-)=-A_- e^{-2zA_+},
\ee
\be
  [A,A_+ ]= \frac{e^{2 z A_+} -1  } z,\quad
[A,A_-]=-2 A_- +z A^2,\quad  [A_+  ,A_- ]= A.
\label{ah}
\ee
The quantum Casimir belonging to the centre of  $U_z sl(2,\R)$ is
now:
\be
{\cal  C}_z=\frac 12 A\, e^{-2zA_+}A+\frac {1-e^{-2zA_+}}{2z}A_-
+A_-\frac {1-e^{-2zA_+}}{2z} + e^{-2zA_+}-1.
\label{aah}
\ee
Obviously, both Lie bialgebras generated by
$\{J_3,J_+,J_-\}$ and
by $\{A,A_+,A_-\}$ are formally
identical, and the classical $r$-matrix  associated to $U_z
sl(2,\R)$ is also
$r= z\,A\wedge A_+$.

Now it is important
to realize that  $\{A,A_+\}$ generates a Hopf subalgebra.
Moreover, this Hopf subalgebra has a known  universal $R$-matrix,
as it was proven in \cite{TT}.

The main result of this letter can be stated as
follows:

\noindent
{\bf Proposition.}
\begin{em} The element
\be
{\cal  R}=\exp\{-z A_+\otimes A\}\exp\{z A\otimes A_+\} .
\label{ai}
\ee
is a quantum universal $R$-matrix for $U_z sl(2,\R)$.
\end{em}

\noindent
{\sl Proof}. The element (\ref{ai}) coincides with  the universal
$R$-matrix for the Hopf subalgebra $\{A,A_+\}$ given in
\cite{TT}. That (\ref{ai}) is a solution of the quantum YBE is a
consequence of the properties of the universal $T$-matrix
\cite{FG} from which ${\cal R}$ was obtained. Thus, we only need
to prove  that property (\ref{aaa}) is fulfilled for the
remaining generator $A_-$. In order to compute the right-hand
side of (\ref{aaa}) we have to take into account  that \be
e^{f}\,\Delta(A_-)\,e^{-f}=\Delta(A_-) +\sum_{n=1}^\infty \frac
1{n!}\,[f,\dots [f,\Delta(A_-)]^{n)}\dots].
\label{aj}
\ee
By setting $f=z\,A\otimes A_+$ we find
\bea
&&\!\!\!\!\!\!\!\!
\!\!\!\![z A\otimes A_+,\Delta(A_-)]= z A\otimes A
+z(-2 A_- +z A^2)\otimes
A_+e^{2 z A_+},\cr
&&\!\!\!\!\!\!\!\!\!\!\!\![z A\otimes A_+,
[z A\otimes A_+,\Delta(A_-)]]=
z A^2\otimes (1- e^{2 z A_+}) - 2 z^2(-2 A_- +z A^2)\otimes
A_+^2e^{2 z A_+}.\cr
&&\label{ak}
\eea
It is now easy to check that, for $n\ge 3$,
\be
[z A\otimes A_+,\dots
[z A\otimes A_+,\Delta(A_-)]^{n)}\dots]
=(-2)^{n-1}z^n (-2 A_- +z A^2)\otimes
A_+^ne^{2 z A_+}.
\label{akb}
\ee
Therefore, if we call $g$ to the following expression
\bea
&&\!\!\!\!\exp\{z A\otimes A_+\}
 \Delta(A_-)\exp\{-z A\otimes A_+\} =
 \Delta(A_-) +  z A\otimes A
+\frac z2 A^2\otimes (1- e^{2 z A_+}) \cr
&&\qquad\qquad -\frac 12 \sum_{n=1}^\infty \frac {(-2z)^n}{n!}
(-2 A_- +z A^2)\otimes
A_+^ne^{2 z A_+}\cr
&& = 1 \otimes A_- + A_-\otimes
e^{2 z A_+ }+z A\otimes A
+\frac z2 A^2\otimes (1- e^{2 z A_+})\cr
&&\qquad\qquad -\frac 12 (-2 A_- +z A^2)\otimes
(e^{-2 z A_+}-1)e^{2 z A_+}\cr
&&=1 \otimes A_- + A_-\otimes 1 + z A\otimes A
=g   ,
\label{al}
\eea
and we compute (\ref{aj}) with $f=-z\,A_+\otimes A$ and $g$
instead of $\Delta(A_-)$, we shall obtain
\bea
&&\!\!\!\!\!\!\!\!
\!\!\!\![-z A_+\otimes A,g ]=-z A_+\otimes(-2 A_- +z A^2)
- z A\otimes A  - z (1- e^{2 z A_+})\otimes A^2,\cr
&&\!\!\!\!\!\!\!\!\!\!\!\!
[-z A_+\otimes A,[-z A_+\otimes A,g ]]=
-2z^2 A_+^2\otimes(-2 A_- +z A^2)+z (1- e^{2 z A_+})
\otimes A^2,\cr
&&\label{am}
\eea
and, for $n\ge 3$,
\be
\!\!\!\!\!\!\!\!\!\!\!\![-z A_+\otimes A,\dots
[-z A_+\otimes A,g ]^{n)}\dots]
=-2^{n-1}z^n A_+^n\otimes
(-2 A_- +z A^2).
\label{amb}
\ee
Now, the proof follows:
\bea
&&\exp\{-z A_+\otimes A\}\, g \,\exp\{z A_+\otimes A\}=g
- z A\otimes A  - \frac z2 (1- e^{2 z A_+})\otimes A^2 \cr
&&\qquad\qquad  -\frac 12 \sum_{n=1}^\infty
\frac {(2zA_+)^n}{n!}
\otimes (-2 A_- +z A^2)\cr
&&=1 \otimes A_- + A_-\otimes 1- \frac z2
(1- e^{2 z A_+})\otimes A^2
  -\frac 12 ( e^{2 z A_+}-1)
\otimes  (-2 A_- +z A^2)\cr
&&=  A_-\otimes 1 + e^{2 z A_+}\otimes A_- =
 \sigma\circ \Delta(A_-).
\label{an}
\eea
Note that ${\cal R}^{-1}={\cal
R}_{21}$ and, therefore, $\cal R$ is a triangular $R$-matrix.

A short disgression concerning some  representations of $U_z
sl(2,\R)$ can be meaningful now. Firstly, the two-dimensional
matrix representation of
$sl(2,\R)$ defined by
\be
D(A)=\left(\begin{array}{ll}
 1 & 0  \\
 0 & -1 \end{array}\right),\quad
D(A_+)=\left(\begin{array}{cc}
 0 & 1  \\
 0 & 0 \end{array}\right),\quad
D(A_-)=\left(\begin{array}{cc}
 0 & 0  \\
 1 & 0 \end{array}\right),\label{ao}
\ee
is also a matrix representation for $U_z sl(2,\R)$ and it can be
used in the FRT approach \cite{FRT} to obtain the quantum group
$Fun_z(SL(2,\R))$. In such a representation, the universal
$R$-matrix (\ref{ai}) has the form:
\bea
&&\!\!\!\!\!\!\!\!
D({\cal R})=I\otimes I + z D(A)\otimes D(A_+)
 -z D(A_+)\otimes D(A) \cr
&&\qquad -z^2 D(A_+) D(A) \otimes D(A) D(A_+)=
\left(\begin{array}{rrrr}
 1 & z & -z& z^2  \\
 0 & 1 & 0& z   \\
 0 & 0 & 1 & -z  \\
 0 & 0 & 0& 1    \end{array}\right).
\label{ap}
\eea
This is exactly the $R$-matrix introduced by  Zakrzewski in
\cite{Zak} and the corresponding quantum $SL(2,\R)$ group is the
one given in \cite{Zak,Ohn}.

On the contrary, the three-dimensional matrix  representation of
$U_z sl(2,\R)$ does not coincide with the classical one and, in
this sense, it can be considered as the first ``non-trivial"
case. By assuming that the deformed natrix realization of the
primitive generator is
\be
D_z(A_+)=\left(\begin{array}{lll}
 0& 1& 0  \\
 0&0& 1 \\
 0&0& 0   \end{array}\right),
\label{aaq}
\ee
and by imposing the commutation rules (\ref{ah}) to be fulfilled, a
straightforward computation leads to the matrices
\be
D_z(A)=\left(\begin{array}{lll}
 2& a& b  \\
 0&0& a - 2\,z \\
 0& 0& -2   \end{array}\right),\quad
D_z(A_-)=\left(\begin{array}{lll}
 - a + 2\,z &  c & d  \\
 2 & 2\,z & b+c \\
 0 & 2 & a  \end{array}\right),\label{aq}
\ee
where  $a$ and $b$ are arbitrary and $c,d$
are given by
\bea
&& c=\frac 14(2\,a\,z - 2\,b -a^2),\cr
&& d=\frac 14(a^2\,z + 2\,b\,z - 2\,a\,z^2 - 2\,a\,b ).
\label{ar}
\eea
Since $D_z(A_+)^3$ vanishes,  the universal
$R$-matrix (\ref{ai}) is realized as
\bea
&&D_z({\cal R})=(1-zD_z(A_+)\otimes D_z(A)+
\frac{z^2}2 D_z(A_+)^2\otimes
D_z(A)^2)\cr &&\qquad\times
(1+zD_z(A)\otimes D_z(A_+)+\frac{z^2}2 D_z(A)^2
\otimes D_z(A_+)^2).
\eea
Explicitly, $D_z({\cal R})$ is the following solution of the
quantum YBE  {\footnotesize{
\be
=\left(\begin{array}{lllllllll}
1 & 2\,z & 2\,{z^2} & -2\,z & 0 & -  b\,z    + a\,{z^2}
 & 2\,{z^2} & b\,z - a\,{z^2} & 0 \cr 0 & 1 & 2\,z &
0 & 0 & 2\,{z^2} & 0
 & 0 & b\,z - a\,{z^2} + 2\,{z^3} \cr 0 & 0 & 1 & 0 &
0 & 2\,z & 0 & 0 &
 2\,{z^2} \cr 0 & 0 & 0 & 1 & 0 & 0 & -2\,z & 2\,{z^2} &
 -  b\,z   + a\,{z^2} - 2\,{z^3} \cr 0 & 0 & 0 & 0 & 1 & 0 & 0
 & 0 & 0 \cr 0 & 0 & 0 & 0 & 0 & 1 & 0 & 0 & 2\,z \cr
0 & 0 & 0 & 0 & 0
 & 0 & 1 & -2\,z & 2\,{z^2} \cr 0 & 0 & 0 & 0 & 0 & 0
& 0 & 1 & -2\,z
 \cr 0 & 0 & 0 & 0 & 0 & 0 & 0 & 0 & 1
\end{array}\right).
\ee}}
Actually, if we define $p=-b\,z+a\,z^2$, we can reduce this
solution to a two parameter $R$-matrix.

Note also that a differential realization of the
commutation rules (\ref{ah}) (with $\lambda\,(\lambda +
1/2)$ as the induced eigenvalue of the quantum
Casimir (\ref{aah})) can be obtained by defining
\bea
&& A=\frac{e^{2 z x} -1  }{z}\,\partial_x -
\lambda\,\frac{e^{2 z x} +1  }{2},
\qquad A_+=x,  \cr
&& A_-=-\frac{e^{2 z x} -1  }{2\,z}\,\partial_x^2 +
\lambda\,\frac{e^{2 z x} +1  }{2}\,\partial_x - z \,
\lambda^2\,\frac{e^{2 z x}
-1  }{8}.
\label{az}
\eea
The $z\to 0$ limit of (\ref{az}) is the usual second
order differential realization of $sl(2,\R)$
\be
A=2\, x\,\partial_x - \lambda,\quad
A_+=x,\label{azz},\quad
A_-=- x\,\partial_x^2 + \lambda\,\partial_x ,
\ee
with the same eigenvalue $\lambda\,(\lambda + 1/2)$ coming from
the classical limit of (\ref{aah}).

Finally, we would like to stress again the fact that the use of
different basis of the same quantum algebra can be helpful in
order to find the associated quantum $R$-matrix. On the other
hand, a complete and systematic study of the representation
theory of this non-standard deformation seems worth to be done in
order to find new matrix solutions of the quantum YBE linked to
$U_z sl(2,\R)$.


\bigskip

This work has been
partially supported by DGICYT (Project PB94-1115) from the
Ministerio de Educaci\'on y Ciencia de Espa\~na.


\bigskip

\noindent
Note added in proof:

After submission of this paper we have been informed of the
paper \cite{iran}, where another expression for a quantum
$R$-matrix of the non-standard deformation of $sl(2,\R)$ has
been given.

\bigskip




\begin{thebibliography}{40}




\bibitem{Yang}  Yang C N 1967
{\it Phys. Rev. Lett.} {\bf 19} 1312

\bibitem{Baxter} Baxter R J  1972
{\it Ann. Phys.} {\bf 70} 193

\bibitem{Gomez}  \'Alvarez--Gaum\'e L, G\'omez C and
Sierra G 1989
{\it Phys. Lett.} {\bf B220} 142

\bibitem{Turaev}  Reshetikhin N Yu and Turaev V G 1990
{\it Commun. Math. Phys.} {\bf 127} 1

\bibitem {FRT}
  Reshetikhin N Y,  Takhtadzhyan L A and    Faddeev L D 1990
{\it Leningrad Math. J.} {\bf 1}  193

\bibitem{Dr} Drinfeld V G 1986
Proceedings of the International Congress of Mathematics,
MRSI Berkeley, 798

\bibitem{Kirillov} Kirillov A N and Reshetikhin N Yu 1988
{\it Preprint LOMI-E-9-88}

\bibitem{Ohn} Ohn C  1992
{\it Lett. Math. Phys.} {\bf 25} 85

\bibitem{Demidov} Demidov E E, Manin Yu I, Mukhin E E and
Zhdanovich D V 1990
{\it Progr. Theor. Phys. Suppl.} {\bf 102} 203

\bibitem{AKS} Aghamohammadi A, Khorrami M and Shariati A  1995
{\it J. Phys. A: Math. Gen.} {\bf 28} L225

\bibitem{Beyond}
 Ballesteros A,  Herranz F J,  del Olmo M A and
Santander M  1995
{\it J. Phys. A: Math. Gen.} {\bf 28} 941

\bibitem{AA}
Vladimirov A A 1993 {\it Mod. Phys. Lett. A} {\bf 8}   2573

\bibitem{RR}
 Ballesteros A, Celeghini E, Herranz F J,  del Olmo M A and
Santander M  1995
{\it J. Phys. A: Math. Gen.}  {\bf 28} 3129

\bibitem{TT}
 Ballesteros A,  Herranz F J,  del Olmo M A, Pere\~na C M and
Santander M
1995 {\it J. Phys. A: Math. Gen.}   {\bf 28} 7113

\bibitem{FG} Fronsdal C and  Galindo A 1993
{\it Lett. Math. Phys.} {\bf 27} 39

\bibitem{Zak} Zakrzewski S 1991
{\it Lett. Math. Phys.} {\bf 22} 287


\bibitem{iran}
  Shariati A, Aghamohammadi A and Khorrami M 1996
 {\it Mod. Phys. Lett. A} {\bf 11} 187


\end{thebibliography}
\end{document}